\documentclass{amsart}
\usepackage{amsfonts}

\newtheorem{theorem}{Theorem}
\theoremstyle{plain}

\newtheorem{corollary}{Corollary}

\newtheorem{lemma}{Lemma}

\newtheorem{proposition}{Proposition}
\newtheorem{remark}{Remark}

\numberwithin{equation}{section}
\input{tcilatex}

\begin{document}
\title[Bessel and Gr\"{u}ss Inequalities for Orthornormal Families]{On
Bessel and Gr\"{u}ss Inequalities for Orthornormal Families in Inner Product
Spaces}
\author{S.S. Dragomir}
\address{School of Computer Science and Mathematics\\
Victoria University of Technology\\
PO Box 14428, MCMC 8001\\
Victoria, Australia.}
\date{May 19, 2003}
\email{sever.dragomir@vu.edu.au}
\urladdr{http://rgmia.vu.edu.au/SSDragomirWeb.html}

\begin{abstract}
A new counterpart of Bessel's inequality for orthornormal families in real
or complex inner product spaces is obtained. Applications for some Gr\"{u}ss
type results are also provided.
\end{abstract}

\maketitle

\section{Introduction}

In the recent paper \cite{SSD2}, the following refinement of the Gr\"{u}ss
inequality was proved.

\begin{theorem}
\label{t1.1}Let $\left( H,\left\langle \cdot ,\cdot \right\rangle \right) $
be an inner product space over $\mathbb{K}$ $\left( \mathbb{K}=\mathbb{R},%
\mathbb{C}\right) $ and $e\in H,$ $\left\Vert e\right\Vert =1.$ If $\phi
,\Phi ,\gamma ,\Gamma $ are real or complex numbers and $x,y$ are vectors in 
$H$ such that either%
\begin{equation}
\func{Re}\left\langle \Phi e-x,x-\phi e\right\rangle \geq 0\text{ \ and \ }%
\func{Re}\left\langle \Gamma e-y,y-\gamma e\right\rangle \geq 0  \label{1.1}
\end{equation}%
or, equivalently, 
\begin{equation}
\left\Vert x-\frac{\phi +\Phi }{2}e\right\Vert \leq \frac{1}{2}\left\vert
\Phi -\phi \right\vert ,\ \ \ \left\Vert y-\frac{\gamma +\Gamma }{2}%
e\right\Vert \leq \frac{1}{2}\left\vert \Gamma -\gamma \right\vert ,
\label{1.2}
\end{equation}%
hold, then we have the following refinement of the Gr\"{u}ss inequality%
\begin{align}
& \left\vert \left\langle x,y\right\rangle -\left\langle x,e\right\rangle
\left\langle e,y\right\rangle \right\vert  \label{1.3} \\
& \leq \frac{1}{4}\left\vert \Phi -\phi \right\vert \left\vert \Gamma
-\gamma \right\vert -\left[ \func{Re}\left\langle \Phi e-x,x-\phi
e\right\rangle \right] ^{\frac{1}{2}}\left[ \func{Re}\left\langle \Gamma
e-y,y-\gamma e\right\rangle \right] ^{\frac{1}{2}}  \notag \\
& \leq \frac{1}{4}\left\vert \Phi -\phi \right\vert \left\vert \Gamma
-\gamma \right\vert .  \notag
\end{align}%
The constant $\frac{1}{4}$ is best possible in both inequalities.
\end{theorem}

Note that the inequality between the first and last term in (\ref{1.3}) was
firstly established in \cite{SSD1}.

A generalisation of the above result for finite families of orthornormal
vectors has been pointed out in \cite{SSD3}.

\begin{theorem}
\label{t1.2}Let $\left\{ e_{i}\right\} _{i\in I}$ be a family of
orthornormal vectors in $H,$ $F$ a finite part of $I,$ $\phi _{i},\Phi
_{i},\gamma _{i},\Gamma _{i}\in \mathbb{K}$ $\left( \mathbb{K}=\mathbb{R},%
\mathbb{C}\right) ,$ $i\in F$ and $x,y\in H.$ If%
\begin{align}
\func{Re}\left\langle \sum_{i=1}^{n}\Phi _{i}e_{i}-x,x-\sum_{i=1}^{n}\phi
_{i}e_{i}\right\rangle & \geq 0,\   \label{1.4} \\
\func{Re}\left\langle \sum_{i=1}^{n}\Gamma
_{i}e_{i}-y,y-\sum_{i=1}^{n}\gamma _{i}e_{i}\right\rangle & \geq 0,  \notag
\end{align}%
or, equivalently,%
\begin{align}
\left\Vert x-\sum_{i\in F}\frac{\Phi _{i}+\phi _{i}}{2}e_{i}\right\Vert &
\leq \frac{1}{2}\left( \sum_{i\in F}\left\vert \Phi _{i}-\phi
_{i}\right\vert ^{2}\right) ^{\frac{1}{2}},  \label{1.5} \\
\left\Vert y-\sum_{i\in F}\frac{\Gamma _{i}+\gamma _{i}}{2}e_{i}\right\Vert
& \leq \frac{1}{2}\left( \sum_{i\in F}\left\vert \Gamma _{i}-\gamma
_{i}\right\vert ^{2}\right) ^{\frac{1}{2}},  \notag
\end{align}%
hold, then we have the inequalities%
\begin{eqnarray}
&&\left\vert \left\langle x,y\right\rangle -\sum_{i\in F}\left\langle
x,e_{i}\right\rangle \left\langle e_{i},y\right\rangle \right\vert
\label{1.6} \\
&\leq &\frac{1}{4}\left( \sum_{i\in F}\left\vert \Phi _{i}-\phi
_{i}\right\vert ^{2}\right) ^{\frac{1}{2}}\cdot \left( \sum_{i\in
F}\left\vert \Gamma _{i}-\gamma _{i}\right\vert ^{2}\right) ^{\frac{1}{2}} 
\notag \\
&&-\left[ \func{Re}\left\langle \sum_{i\in F}\Phi _{i}e_{i}-x,x-\sum_{i\in
F}\phi _{i}e_{i}\right\rangle \right] ^{\frac{1}{2}}\left[ \func{Re}%
\left\langle \sum_{i\in F}\Gamma _{i}e_{i}-y,y-\sum_{i\in F}\gamma
_{i}e_{i}\right\rangle \right] ^{\frac{1}{2}}  \notag \\
&\leq &\frac{1}{4}\left( \sum_{i\in F}\left\vert \Phi _{i}-\phi
_{i}\right\vert ^{2}\right) ^{\frac{1}{2}}\cdot \left( \sum_{i\in
F}\left\vert \Gamma _{i}-\gamma _{i}\right\vert ^{2}\right) ^{\frac{1}{2}}. 
\notag
\end{eqnarray}%
The constant $\frac{1}{4}$ is best possible.
\end{theorem}

\begin{remark}
We note that the inequality between the first term and the last term for
real inner products under the assumption (\ref{1.4}), has been proved by N.
Ujevi\'{c} in \cite{NU}.
\end{remark}

The following corollary provides a counterpart for the well known Bessel's
inequality in real or complex inner product spaces (see also \cite{SSD3}).

\begin{corollary}
\label{c1.3}With the above assumptions for $\left\{ e_{i}\right\} _{i\in
I},F,\phi _{i},$ $\Phi _{i}$ and $x,$ one has the inequalities%
\begin{align}
0& \leq \left\Vert x\right\Vert ^{2}-\sum_{i\in F}\left\vert \left\langle
x,e_{i}\right\rangle \right\vert ^{2}  \label{1.7} \\
& \leq \frac{1}{4}\sum_{i\in F}\left\vert \Phi _{i}-\phi _{i}\right\vert
^{2}-\func{Re}\left\langle \sum_{i\in F}\Phi _{i}e_{i}-x,x-\sum_{i\in F}\phi
_{i}e_{i}\right\rangle  \notag \\
& \leq \frac{1}{4}\sum_{i\in F}\left\vert \Phi _{i}-\phi _{i}\right\vert
^{2},  \notag
\end{align}%
with $\frac{1}{4}$ as the best possible constant in both inequalities.
\end{corollary}

It is the main aim of this paper to point out a different counterpart for
the Bessel and Gr\"{u}ss inequalities stated above. Some related results are
also outlined.

\section{A Counterpart of Bessel's Inequality}

The following lemma holds.

\begin{lemma}
\label{l2.1}Let $\left\{ e_{i}\right\} _{i\in I}$ be a family of
orthornormal vectors in $H,$ $F$ a finite part of $I,$ $\lambda _{i}\in 
\mathbb{K},$ $i\in F$, $r>0$ and $x\in H.$ If%
\begin{equation}
\left\Vert x-\sum_{i\in F}\lambda _{i}e_{i}\right\Vert \leq r,  \label{2.1}
\end{equation}%
then we have the inequality%
\begin{equation}
0\leq \left\Vert x\right\Vert ^{2}-\sum_{i\in F}\left\vert \left\langle
x,e_{i}\right\rangle \right\vert ^{2}\leq r^{2}-\sum_{i\in F}\left\vert
\lambda _{i}-\left\langle x,e_{i}\right\rangle \right\vert ^{2}.  \label{2.2}
\end{equation}
\end{lemma}

\begin{proof}
Consider%
\begin{align*}
I_{1}& :=\left\Vert x-\sum_{i\in F}\lambda _{i}e_{i}\right\Vert
^{2}=\left\langle x-\sum_{i\in F}\lambda _{i}e_{i},x-\sum_{j\in F}\lambda
_{j}e_{j}\right\rangle \\
& =\left\Vert x\right\Vert ^{2}-\sum_{i\in F}\lambda _{i}\overline{%
\left\langle x,e_{i}\right\rangle }-\sum_{i\in F}\overline{\lambda _{i}}%
\left\langle x,e_{i}\right\rangle +\sum_{i\in F}\sum_{j\in F}\lambda _{i}%
\overline{\lambda _{j}}\left\langle e_{i},e_{j}\right\rangle \\
& =\left\Vert x\right\Vert ^{2}-\sum_{i\in F}\lambda _{i}\overline{%
\left\langle x,e_{i}\right\rangle }-\sum_{i\in F}\overline{\lambda _{i}}%
\left\langle x,e_{i}\right\rangle +\sum_{i\in F}\left\vert \lambda
_{i}\right\vert ^{2}
\end{align*}%
and%
\begin{align*}
I_{2}& :=\sum_{i\in F}\left\vert \lambda _{i}-\left\langle
x,e_{i}\right\rangle \right\vert ^{2}=\sum_{i\in F}\left( \lambda
_{i}-\left\langle x,e_{i}\right\rangle \right) \left( \overline{\lambda _{i}}%
-\overline{\left\langle x,e_{i}\right\rangle }\right) \\
& =\sum_{i\in F}\left[ \left\vert \lambda _{i}\right\vert ^{2}+\left\vert
\left\langle x,e_{i}\right\rangle \right\vert ^{2}-\overline{\lambda _{i}}%
\left\langle x,e_{i}\right\rangle -\lambda _{i}\overline{\left\langle
x,e_{i}\right\rangle }\right] \\
& =\sum_{i\in F}\left\vert \lambda _{i}\right\vert ^{2}+\sum_{i\in
F}\left\vert \left\langle x,e_{i}\right\rangle \right\vert ^{2}-\sum_{i\in F}%
\overline{\lambda _{i}}\left\langle x,e_{i}\right\rangle -\sum_{i\in
F}\lambda _{i}\overline{\left\langle x,e_{i}\right\rangle }.
\end{align*}%
If we subtract $I_{2}$ from $I_{1}$ we deduce an equality that is
interesting in its own right%
\begin{equation}
\left\Vert x-\sum_{i\in F}\lambda _{i}e_{i}\right\Vert ^{2}-\sum_{i\in
F}\left\vert \lambda _{i}-\left\langle x,e_{i}\right\rangle \right\vert
^{2}=\left\Vert x\right\Vert ^{2}-\sum_{i\in F}\left\vert \left\langle
x,e_{i}\right\rangle \right\vert ^{2},  \label{2.3}
\end{equation}%
from which we easily deduce (\ref{2.2}).
\end{proof}

The following counterpart of Bessel's inequality holds.

\begin{theorem}
\label{t2.2}Let $\left\{ e_{i}\right\} _{i\in I}$ be a family of
orthornormal vectors in $H,$ $F$ a finite part of $I,$ $\phi _{i},$ $\Phi
_{i},$ $i\in I$ real or complex numbers. For a $x\in H,$ if either

\begin{enumerate}
\item[(i)] $\func{Re}\left\langle \sum_{i\in F}\Phi _{i}e_{i}-x,x-\sum_{i\in
F}\phi _{i}e_{i}\right\rangle \geq 0;$\newline

or, equivalently,

\item[(ii)] $\left\Vert x-\sum_{i\in F}\frac{\phi _{i}+\Phi _{i}}{2}%
e_{i}\right\Vert \leq \frac{1}{2}\left( \sum_{i\in F}\left\vert \Phi
_{i}-\phi _{i}\right\vert ^{2}\right) ^{\frac{1}{2}};$
\end{enumerate}

holds, then the following counterpart of Bessel's inequality%
\begin{align}
0& \leq \left\Vert x\right\Vert ^{2}-\sum_{i\in F}\left\vert \left\langle
x,e_{i}\right\rangle \right\vert ^{2}  \label{2.4} \\
& \leq \frac{1}{4}\sum_{i\in F}\left\vert \Phi _{i}-\phi _{i}\right\vert
^{2}-\sum_{i\in F}\left\vert \frac{\phi _{i}+\Phi _{i}}{2}-\left\langle
x,e_{i}\right\rangle \right\vert ^{2}  \notag \\
& \leq \frac{1}{4}\sum_{i\in F}\left\vert \Phi _{i}-\phi _{i}\right\vert
^{2},  \notag
\end{align}%
is valid.

The constant $\frac{1}{4}$ is best possible in both inequalities.
\end{theorem}

\begin{proof}
Firstly, we observe that for $y,a,A\in H,$ the following are equivalent%
\begin{equation}
\func{Re}\left\langle A-y,y-a\right\rangle \geq 0  \label{2.5}
\end{equation}%
and%
\begin{equation}
\left\Vert y-\frac{a+A}{2}\right\Vert \leq \frac{1}{2}\left\Vert
A-a\right\Vert .  \label{2.6}
\end{equation}%
Now, for $a=\sum_{i\in F}\phi _{i}e_{i},$ $A=\sum_{i\in F}\Phi _{i}e_{i},$
we have%
\begin{align*}
\left\Vert A-a\right\Vert & =\left\Vert \sum_{i\in F}\left( \Phi _{i}-\phi
_{i}\right) e_{i}\right\Vert =\left( \left\Vert \sum_{i\in F}\left( \Phi
_{i}-\phi _{i}\right) e_{i}\right\Vert ^{2}\right) ^{\frac{1}{2}} \\
& =\left( \sum_{i\in F}\left\vert \Phi _{i}-\phi _{i}\right\vert
^{2}\left\Vert e_{i}\right\Vert ^{2}\right) ^{\frac{1}{2}}=\left( \sum_{i\in
F}\left\vert \Phi _{i}-\phi _{i}\right\vert ^{2}\right) ^{\frac{1}{2}},
\end{align*}%
giving, for $y=x,$ the desired equivalence.

Now, if we apply Lemma \ref{l2.1} for $\lambda _{i}=\frac{\phi _{i}+\Phi _{i}%
}{2}$ and 
\begin{equation*}
r:=\frac{1}{2}\left( \sum_{i\in F}\left\vert \Phi _{i}-\phi _{i}\right\vert
^{2}\right) ^{\frac{1}{2}},
\end{equation*}%
we deduce the first inequality in (\ref{2.4}).

Let us prove that $\frac{1}{4}$ is best possible in the second inequality in
(\ref{2.4}). Assume that there is a $c>0$ such that%
\begin{equation}
0\leq \left\Vert x\right\Vert ^{2}-\sum_{i\in F}\left\vert \left\langle
x,e_{i}\right\rangle \right\vert ^{2}\leq c\sum_{i\in F}\left\vert \Phi
_{i}-\phi _{i}\right\vert ^{2}-\sum_{i\in F}\left\vert \frac{\phi _{i}+\Phi
_{i}}{2}-\left\langle x,e_{i}\right\rangle \right\vert ^{2},  \label{2.7}
\end{equation}%
provided that $\phi _{i},$ $\Phi _{i},x$ and $F$ satisfy (i) and (ii).

We choose $F=\left\{ 1\right\} ,$ $e_{1}=e=\left( \frac{1}{\sqrt{2}},\frac{1%
}{\sqrt{2}}\right) \in \mathbb{R}^{2},$ $x=\left( x_{1},x_{2}\right) \in 
\mathbb{R}^{2},$ $\Phi _{1}=\Phi =m>0,$ $\phi _{1}=\phi =-m,$ $H=\mathbb{R}%
^{2}$ to get from (\ref{2.7}) that%
\begin{align}
0& \leq x_{1}^{2}+x_{2}^{2}-\frac{\left( x_{1}+x_{2}\right) ^{2}}{2}
\label{2.8} \\
& \leq 4cm^{2}-\frac{\left( x_{1}+x_{2}\right) ^{2}}{2},  \notag
\end{align}%
provided%
\begin{align}
0& \leq \left\langle me-x,x+me\right\rangle  \label{2.9} \\
& =\left( \frac{m}{\sqrt{2}}-x_{1}\right) \left( x_{1}+\frac{m}{\sqrt{2}}%
\right) +\left( \frac{m}{\sqrt{2}}-x_{2}\right) \left( x_{2}+\frac{m}{\sqrt{2%
}}\right) .  \notag
\end{align}%
From (\ref{2.8}) we get%
\begin{equation}
x_{1}^{2}+x_{2}^{2}\leq 4cm^{2}  \label{2.10}
\end{equation}%
provided (\ref{2.9}) holds.

If we choose $x_{1}=\frac{m}{\sqrt{2}},$ $x_{2}=-\frac{m}{\sqrt{2}},$ then (%
\ref{2.9}) is fulfilled and by (\ref{2.10}) we get $m^{2}\leq 4cm^{2},$
giving $c\geq \frac{1}{4}.$
\end{proof}

\begin{remark}
\label{r2.3}If $F=\left\{ 1\right\} ,$ $e_{1}=1,$ $\left\Vert e\right\Vert =1
$ and for $\phi ,\Phi \in \mathbb{K}$ and $x\in H$ one has either%
\begin{equation}
\func{Re}\left\langle \Phi e-x,x-\phi e\right\rangle \geq 0  \label{2.11}
\end{equation}%
or, equivalently,%
\begin{equation}
\left\Vert x-\frac{\phi +\Phi }{2}e\right\Vert \leq \frac{1}{2}\left\vert
\Phi -\phi \right\vert ,  \label{2.12}
\end{equation}%
then%
\begin{align}
0& \leq \left\Vert x\right\Vert ^{2}-\left\vert \left\langle
x,e\right\rangle \right\vert ^{2}  \label{2.13} \\
& \leq \frac{1}{4}\left\vert \Phi -\phi \right\vert ^{2}-\left\vert \frac{%
\phi +\Phi }{2}-\left\langle x,e\right\rangle \right\vert ^{2}\leq \frac{1}{4%
}\left\vert \Phi -\phi \right\vert ^{2}.  \notag
\end{align}%
The constant $\frac{1}{4}$ is best possible in both inequalities.
\end{remark}

\begin{remark}
\label{r2.4}It is important to compare the bounds provided by Corollary \ref%
{c1.3} and Theorem \ref{t2.2}.

For this purpose, consider%
\begin{equation*}
B_{1}\left( x,e,\phi ,\Phi \right) :=\frac{1}{4}\left( \Phi -\phi \right)
^{2}-\left\langle \Phi e-x,x-\phi e\right\rangle
\end{equation*}%
and%
\begin{equation*}
B_{2}\left( x,e,\phi ,\Phi \right) :=\frac{1}{4}\left( \Phi -\phi \right)
^{2}-\left( \frac{\phi +\Phi }{2}-\left\langle x,e\right\rangle \right) ^{2},
\end{equation*}%
where $H$ is a real inner product, $e\in H,$ $\left\Vert e\right\Vert =1,$ $%
x\in H,$ $\phi ,\Phi \in \mathbb{R}$ with 
\begin{equation*}
\left\langle \Phi e-x,x-\phi e\right\rangle \geq 0,
\end{equation*}%
or equivalently, 
\begin{equation*}
\left\Vert x-\frac{\phi +\Phi }{2}e\right\Vert \leq \frac{1}{2}\left\vert
\Phi -\phi \right\vert .
\end{equation*}

If we choose $\phi =-1,$ $\Phi =1,$ then we have%
\begin{align*}
B_{1}\left( x,e\right) & =1-\left\langle e-x,x+e\right\rangle =1-\left(
\left\Vert e\right\Vert ^{2}-\left\Vert x\right\Vert ^{2}\right) =\left\Vert
x\right\Vert ^{2}, \\
B_{2}\left( x,e\right) & =1-\left\langle x,e\right\rangle ^{2},
\end{align*}%
provided $\left\Vert x\right\Vert \leq 1.$

Choose $x=ke,$ with $0<k\leq 1.$ Then we get%
\begin{equation*}
B_{1}\left( k\right) =k^{2},\ \ B_{2}\left( k\right) =1-k^{2},
\end{equation*}%
which shows that $B_{1}\left( k\right) >B_{2}\left( k\right) $ if $0<k<\frac{%
\sqrt{2}}{2}$ and $B_{1}\left( k\right) <B_{2}\left( k\right) $ if $\frac{%
\sqrt{2}}{2}<k\leq 1.$
\end{remark}

We may state the following proposition.

\begin{proposition}
\label{p2.5}Let $\left\{ e_{i}\right\} _{i\in I}$ be a family of
orthornormal vectors in $H,$ $F$ a finite part of $I,$ $\phi _{i},$ $\Phi
_{i}\in \mathbb{K}$ $\left( i\in F\right) $. If $x\in H$ either satisfies
(i), or, equivalently, (ii) of Theorem \ref{t2.2}, then the upper bounds%
\begin{align*}
B_{1}\left( x,e,\mathbf{\phi },\mathbf{\Phi },F\right) & :=\frac{1}{4}%
\sum_{i\in F}\left\vert \Phi _{i}-\phi _{i}\right\vert ^{2}-\func{Re}%
\left\langle \sum_{i\in F}\Phi _{i}e_{i}-x,x-\sum_{i\in F}\phi
_{i}e_{i}\right\rangle , \\
B_{2}\left( x,e,\mathbf{\phi },\mathbf{\Phi },F\right) & :=\frac{1}{4}%
\sum_{i\in F}\left\vert \Phi _{i}-\phi _{i}\right\vert ^{2}-\sum_{i\in
F}\left\vert \frac{\phi _{i}+\Phi _{i}}{2}-\left\langle x,e_{i}\right\rangle
\right\vert ^{2},
\end{align*}%
for the Bessel's difference $B_{s}\left( x,\mathbf{e},F\right) :=\left\Vert
x\right\Vert ^{2}-\sum_{i\in F}\left\vert \left\langle x,e_{i}\right\rangle
\right\vert ^{2},$ cannot be compared in general.
\end{proposition}

\section{A Refinement of the Gr\"{u}ss Inequality\label{s3}}

The following result holds.

\begin{theorem}
\label{t3.1}Let $\left\{ e_{i}\right\} _{i\in I}$ be a family of
orthornormal vectors in $H,$ $F$ a finite part of $I$, $\phi _{i},\Phi _{i},$
$\gamma _{i},\Gamma _{i}\in \mathbb{K},\ i\in F$ and $x,y\in H.$ If either%
\begin{align}
\func{Re}\left\langle \sum_{i\in F}\Phi _{i}e_{i}-x,x-\sum_{i\in F}\phi
_{i}e_{i}\right\rangle & \geq 0,  \label{3.1} \\
\func{Re}\left\langle \sum_{i\in F}\Gamma _{i}e_{i}-y,y-\sum_{i\in F}\gamma
_{i}e_{i}\right\rangle & \geq 0,  \notag
\end{align}%
or, equivalently,%
\begin{align}
\left\Vert x-\sum_{i\in F}\frac{\Phi _{i}+\phi _{i}}{2}e_{i}\right\Vert &
\leq \frac{1}{2}\left( \sum_{i\in F}\left\vert \Phi _{i}-\phi
_{i}\right\vert ^{2}\right) ^{\frac{1}{2}},  \label{3.2} \\
\left\Vert y-\sum_{i\in F}\frac{\Gamma _{i}+\gamma _{i}}{2}e_{i}\right\Vert
& \leq \frac{1}{2}\left( \sum_{i\in F}\left\vert \Gamma _{i}-\gamma
_{i}\right\vert ^{2}\right) ^{\frac{1}{2}},  \notag
\end{align}%
hold, then we have the inequalities%
\begin{align}
0& \leq \left\vert \left\langle x,y\right\rangle -\sum_{i\in F}\left\langle
x,e_{i}\right\rangle \left\langle e_{i},y\right\rangle \right\vert
\label{3.3} \\
& \leq \frac{1}{4}\left( \sum_{i\in F}\left\vert \Phi _{i}-\phi
_{i}\right\vert ^{2}\right) ^{\frac{1}{2}}\cdot \left( \sum_{i\in
F}\left\vert \Gamma _{i}-\gamma _{i}\right\vert ^{2}\right) ^{\frac{1}{2}} 
\notag \\
& \ \ \ \ \ \ \ \ \ \ \ \ \ \ \ \ \ -\sum_{i\in F}\left\vert \frac{\Phi
_{i}+\phi _{i}}{2}-\left\langle x,e_{i}\right\rangle \right\vert \left\vert 
\frac{\Gamma _{i}+\gamma _{i}}{2}-\left\langle y,e_{i}\right\rangle
\right\vert  \notag \\
& \leq \frac{1}{4}\left( \sum_{i\in F}\left\vert \Phi _{i}-\phi
_{i}\right\vert ^{2}\right) ^{\frac{1}{2}}\cdot \left( \sum_{i\in
F}\left\vert \Gamma _{i}-\gamma _{i}\right\vert ^{2}\right) ^{\frac{1}{2}}. 
\notag
\end{align}%
The constant $\frac{1}{4}$ is best possible.
\end{theorem}

\begin{proof}
Using Schwartz's inequality in the inner product space $\left(
H,\left\langle \cdot ,\cdot \right\rangle \right) $ one has%
\begin{multline}
\left\vert \left\langle x-\sum_{i\in F}\left\langle x,e_{i}\right\rangle
e_{i},y-\sum_{i\in F}\left\langle y,e_{i}\right\rangle e_{i}\right\rangle
\right\vert ^{2}  \label{3.4} \\
\leq \left\Vert x-\sum_{i\in F}\left\langle x,e_{i}\right\rangle
e_{i}\right\Vert ^{2}\left\Vert y-\sum_{i\in F}\left\langle
y,e_{i}\right\rangle e_{i}\right\Vert ^{2}
\end{multline}%
and since a simple calculation shows that 
\begin{equation*}
\left\langle x-\sum_{i\in F}\left\langle x,e_{i}\right\rangle
e_{i},y-\sum_{i\in F}\left\langle y,e_{i}\right\rangle e_{i}\right\rangle
=\left\langle x,y\right\rangle -\sum_{i\in F}\left\langle
x,e_{i}\right\rangle \left\langle e_{i},y\right\rangle 
\end{equation*}%
and 
\begin{equation*}
\left\Vert x-\sum_{i\in F}\left\langle x,e_{i}\right\rangle e_{i}\right\Vert
^{2}=\left\Vert x\right\Vert ^{2}-\sum_{i\in F}\left\vert \left\langle
x,e_{i}\right\rangle \right\vert ^{2}
\end{equation*}%
for any $x,y\in H,$ then by (\ref{3.4}) and by the counterpart of Bessel's
inequality in Theorem \ref{t2.2}, we have%
\begin{align}
& \left\vert \left\langle x,y\right\rangle -\sum_{i\in F}\left\langle
x,e_{i}\right\rangle \left\langle e_{i},y\right\rangle \right\vert ^{2}
\label{3.5} \\
& \leq \left( \left\Vert x\right\Vert ^{2}-\sum_{i\in F}\left\vert
\left\langle x,e_{i}\right\rangle \right\vert ^{2}\right) \left( \left\Vert
y\right\Vert ^{2}-\sum_{i\in F}\left\vert \left\langle y,e_{i}\right\rangle
\right\vert ^{2}\right)   \notag \\
& \leq \left[ \frac{1}{4}\sum_{i\in F}\left\vert \Phi _{i}-\phi
_{i}\right\vert ^{2}-\sum_{i\in F}\left\vert \frac{\Phi _{i}+\phi _{i}}{2}%
-\left\langle x,e_{i}\right\rangle \right\vert ^{2}\right]   \notag \\
& \ \ \ \ \ \ \ \ \ \ \ \ \ \ \ \ \times \left[ \frac{1}{4}\sum_{i\in
F}\left\vert \Gamma _{i}-\gamma _{i}\right\vert ^{2}-\sum_{i\in F}\left\vert 
\frac{\Gamma _{i}+\gamma _{i}}{2}-\left\langle y,e_{i}\right\rangle
\right\vert ^{2}\right]   \notag \\
& :=K.  \notag
\end{align}%
Using Acz\'{e}l's inequality for real numbers, i.e., we recall that%
\begin{equation}
\left( a^{2}-\sum_{i\in F}a_{i}^{2}\right) \left( b^{2}-\sum_{i\in
F}b_{i}^{2}\right) \leq \left( ab-\sum_{i\in F}a_{i}b_{i}\right) ^{2},
\label{3.6}
\end{equation}%
provided that $a,b,a_{i},b_{i}>0,$ $i\in F,$ we may state that%
\begin{multline}
K\leq \left[ \frac{1}{4}\left( \sum_{i\in F}\left\vert \Phi _{i}-\phi
_{i}\right\vert ^{2}\right) ^{\frac{1}{2}}\cdot \left( \sum_{i\in
F}\left\vert \Gamma _{i}-\gamma _{i}\right\vert ^{2}\right) ^{\frac{1}{2}%
}\right.   \label{3.7} \\
\left. -\sum_{i\in F}\left\vert \frac{\Phi _{i}+\phi _{i}}{2}-\left\langle
x,e_{i}\right\rangle \right\vert \left\vert \frac{\Gamma _{i}+\gamma _{i}}{2}%
-\left\langle y,e_{i}\right\rangle \right\vert \right] ^{2}.
\end{multline}%
Using (\ref{3.5}) and (\ref{3.7}) we conclude that%
\begin{multline}
\left\vert \left\langle x,y\right\rangle -\sum_{i\in F}\left\langle
x,e_{i}\right\rangle \left\langle e_{i},y\right\rangle \right\vert ^{2}\leq 
\left[ \frac{1}{4}\left( \sum_{i\in F}\left\vert \Phi _{i}-\phi
_{i}\right\vert ^{2}\right) ^{\frac{1}{2}}\cdot \left( \sum_{i\in
F}\left\vert \Gamma _{i}-\gamma _{i}\right\vert ^{2}\right) ^{\frac{1}{2}%
}\right.   \label{3.8} \\
-\left. \sum_{i\in F}\left\vert \frac{\Phi _{i}+\phi _{i}}{2}-\left\langle
x,e_{i}\right\rangle \right\vert \left\vert \frac{\Gamma _{i}+\gamma _{i}}{2}%
-\left\langle y,e_{i}\right\rangle \right\vert \right] ^{2}.
\end{multline}%
Taking the square root in (\ref{3.8}) and taking into account that the
quantity in the last square brackets is nonnegative (see for example (\ref%
{2.4})), we deduce the second inequality in (\ref{3.3}).

The fact that $\frac{1}{4}$ is the best possible constant follows by Theorem %
\ref{t2.2} and we omit the details.
\end{proof}

The following corollary may be stated.

\begin{corollary}
\label{c3.2}Let $e\in H,$ $\left\Vert e\right\Vert =1,$ $\phi ,\Phi ,\gamma
,\Gamma \in \mathbb{K}$ and $x,y\in H$ such that either%
\begin{equation}
\func{Re}\left\langle \Phi e-x,x-\phi e\right\rangle \geq 0\text{ \ and \ }%
\func{Re}\left\langle \Gamma e-y,y-\gamma e\right\rangle \geq 0  \label{3.9}
\end{equation}%
or, equivalently,%
\begin{equation}
\left\Vert x-\frac{\phi +\Phi }{2}e\right\Vert \leq \frac{1}{2}\left\vert
\Phi -\phi \right\vert ,\ \ \ \left\Vert y-\frac{\gamma +\Gamma }{2}%
e\right\Vert \leq \frac{1}{2}\left\vert \Gamma -\gamma \right\vert .
\label{3.10}
\end{equation}%
Then we have the following refinement of Gr\"{u}ss' inequality%
\begin{align}
0& \leq \left\vert \left\langle x,y\right\rangle -\left\langle
x,e\right\rangle \left\langle e,y\right\rangle \right\vert   \label{3.11} \\
& \leq \frac{1}{4}\left\vert \Phi -\phi \right\vert \left\vert \Gamma
-\gamma \right\vert -\left\vert \frac{\phi +\Phi }{2}-\left\langle
x,e\right\rangle \right\vert \left\vert \frac{\gamma +\Gamma }{2}%
-\left\langle y,e\right\rangle \right\vert   \notag \\
& \leq \frac{1}{4}\left\vert \Phi -\phi \right\vert \left\vert \Gamma
-\gamma \right\vert .  \notag
\end{align}%
The constant $\frac{1}{4}$ is best possible in both inequalities.
\end{corollary}

\section{Some Companion Inequalities\label{s4}}

The following companion of the Gr\"{u}ss inequality also holds.

\begin{theorem}
\label{t4.1}Let $\left\{ e_{i}\right\} _{i\in I}$ be a family of
orthornormal vectors in $H,$ $F$ a finite part of $I$ and $\phi _{i},\Phi
_{i}\in \mathbb{K},\ i\in F$, $x,y\in H$ and $\lambda \in \left( 0,1\right) ,
$ such that either%
\begin{equation}
\func{Re}\left\langle \sum_{i\in F}\Phi _{i}e_{i}-\left( \lambda x+\left(
1-\lambda \right) y\right) ,\lambda x+\left( 1-\lambda \right) y-\sum_{i\in
F}\phi _{i}e_{i}\right\rangle \geq 0  \label{4.1}
\end{equation}%
or, equivalently,%
\begin{equation}
\left\Vert \lambda x+\left( 1-\lambda \right) y-\sum_{i\in F}\frac{\Phi
_{i}+\phi _{i}}{2}\cdot e_{i}\right\Vert \leq \frac{1}{2}\left( \sum_{i\in
F}\left\vert \Phi _{i}-\phi _{i}\right\vert ^{2}\right) ^{\frac{1}{2}},
\label{4.2}
\end{equation}%
holds. Then we have the inequality%
\begin{align}
& \func{Re}\left[ \left\langle x,y\right\rangle -\sum_{i\in F}\left\langle
x,e_{i}\right\rangle \left\langle e_{i},y\right\rangle \right]   \label{4.3}
\\
& \leq \frac{1}{16}\cdot \frac{1}{\lambda \left( 1-\lambda \right) }%
\sum_{i\in F}\left\vert \Phi _{i}-\phi _{i}\right\vert ^{2}  \notag \\
& \ \ \ \ \ \ \ \ \ \ \ \ \ \ -\frac{1}{4}\frac{1}{\lambda \left( 1-\lambda
\right) }\sum_{i\in F}\left\vert \frac{\Phi _{i}+\phi _{i}}{2}-\left\langle
\lambda x+\left( 1-\lambda \right) y,e_{i}\right\rangle \right\vert ^{2} 
\notag \\
& \leq \frac{1}{16}\cdot \frac{1}{\lambda \left( 1-\lambda \right) }%
\sum_{i\in F}\left\vert \Phi _{i}-\phi _{i}\right\vert ^{2}.  \notag
\end{align}%
The constant $\frac{1}{16}$ is the best possible constant in (\ref{4.3}) in
the sense that it cannot be replaced by a smaller constant.
\end{theorem}

\begin{proof}
We know that for any$\ z,u\in H,$ one has%
\begin{equation*}
\func{Re}\left\langle z,u\right\rangle \leq \frac{1}{4}\left\Vert
z+u\right\Vert ^{2}.
\end{equation*}%
Then for any $a,b\in H$ and $\lambda \in \left( 0,1\right) $ one has%
\begin{equation}
\func{Re}\left\langle a,b\right\rangle \leq \frac{1}{4\lambda \left(
1-\lambda \right) }\left\Vert \lambda a+\left( 1-\lambda \right)
b\right\Vert ^{2}.  \label{4.4}
\end{equation}%
Since%
\begin{equation*}
\left\langle x,y\right\rangle -\sum_{i\in F}\left\langle
x,e_{i}\right\rangle \left\langle e_{i},y\right\rangle =\left\langle
x-\sum_{i\in F}\left\langle x,e_{i}\right\rangle e_{i},y-\sum_{i\in
F}\left\langle y,e_{i}\right\rangle e_{i}\right\rangle ,
\end{equation*}%
for any \thinspace $x,y\in H,$ then, by (\ref{4.4}), we get%
\begin{align}
& \func{Re}\left[ \left\langle x,y\right\rangle -\sum_{i\in F}\left\langle
x,e_{i}\right\rangle \left\langle e_{i},y\right\rangle \right]  \label{4.5}
\\
& =\func{Re}\left[ \left\langle x-\sum_{i\in F}\left\langle
x,e_{i}\right\rangle e_{i},y-\sum_{i\in F}\left\langle y,e_{i}\right\rangle
e_{i}\right\rangle \right]  \notag \\
& \leq \frac{1}{4\lambda \left( 1-\lambda \right) }\left\Vert \lambda \left(
x-\sum_{i\in F}\left\langle x,e_{i}\right\rangle e_{i}\right) +\left(
1-\lambda \right) \left( y-\sum_{i\in F}\left\langle y,e_{i}\right\rangle
e_{i}\right) \right\Vert ^{2}  \notag \\
& =\frac{1}{4\lambda \left( 1-\lambda \right) }\left\Vert \lambda x+\left(
1-\lambda \right) y-\sum_{i\in F}\left\langle \lambda x+\left( 1-\lambda
\right) y,e_{i}\right\rangle e_{i}\right\Vert ^{2}  \notag \\
& =\frac{1}{4\lambda \left( 1-\lambda \right) }\left[ \left\Vert \lambda
x+\left( 1-\lambda \right) y\right\Vert ^{2}-\sum_{i\in F}\left\vert
\left\langle \lambda x+\left( 1-\lambda \right) y,e_{i}\right\rangle
\right\vert ^{2}\right] .  \notag
\end{align}%
If we apply the counterpart of Bessel's inequality in Theorem \ref{t2.2} for 
$\lambda x+\left( 1-\lambda \right) y,$ we may state that%
\begin{align}
& \left\Vert \lambda x+\left( 1-\lambda \right) y\right\Vert ^{2}-\sum_{i\in
F}\left\vert \left\langle \lambda x+\left( 1-\lambda \right)
y,e_{i}\right\rangle \right\vert ^{2}  \label{4.6} \\
& \leq \frac{1}{4}\sum_{i\in F}\left\vert \Phi _{i}-\phi _{i}\right\vert
^{2}-\sum_{i\in F}\left\vert \frac{\Phi _{i}+\phi _{i}}{2}-\left\langle
\lambda x+\left( 1-\lambda \right) y,e_{i}\right\rangle \right\vert ^{2} 
\notag \\
& \leq \frac{1}{4}\sum_{i\in F}\left\vert \Phi _{i}-\phi _{i}\right\vert
^{2}.  \notag
\end{align}%
Now, by making use of (\ref{4.5}) and (\ref{4.6}), we deduce (\ref{4.3}).

The fact that $\frac{1}{16}$ is the best possible constant in (\ref{4.3})
follows by the fact that if in (\ref{4.1}) we choose $x=y,$ then it becomes
(i) of Theorem \ref{t2.2}, implying for $\lambda =\frac{1}{2}$ (\ref{2.4}),
for which, we have shown that $\frac{1}{4}$ was the best constant.
\end{proof}

\begin{remark}
\label{r1}In practical applications we may use only the inequality between
the first and the last term in (\ref{4.3}).
\end{remark}

\begin{remark}
\label{r2}If in Theorem \ref{t4.1}, we choose $\lambda =\frac{1}{2},$ then
we get%
\begin{align}
& \func{Re}\left[ \left\langle x,y\right\rangle -\sum_{i\in F}\left\langle
x,e_{i}\right\rangle \left\langle e_{i},y\right\rangle \right]   \label{4.7}
\\
& \leq \frac{1}{4}\sum_{i\in F}\left\vert \Phi _{i}-\phi _{i}\right\vert
^{2}-\sum_{i\in F}\left\vert \frac{\Phi _{i}+\phi _{i}}{2}-\left\langle 
\frac{x+y}{2},e_{i}\right\rangle \right\vert ^{2}  \notag \\
& \leq \frac{1}{4}\sum_{i\in F}\left\vert \Phi _{i}-\phi _{i}\right\vert
^{2},  \notag
\end{align}%
provided%
\begin{equation*}
\func{Re}\left\langle \sum_{i\in F}\Phi _{i}e_{i}-\frac{x+y}{2},\frac{x+y}{2}%
-\sum_{i\in F}\phi _{i}e_{i}\right\rangle \geq 0
\end{equation*}%
or, equivalently,%
\begin{equation*}
\left\Vert \frac{x+y}{2}-\sum_{i\in F}\frac{\Phi _{i}+\phi _{i}}{2}\cdot
e_{i}\right\Vert \leq \frac{1}{2}\left( \sum_{i\in F}\left\vert \Phi
_{i}-\phi _{i}\right\vert ^{2}\right) ^{\frac{1}{2}}.
\end{equation*}
\end{remark}

We note that (\ref{4.7}) is a refinement of the corresponding result in \cite%
{SSD3}.

\begin{corollary}
\label{c4.2}With the assumptions of Theorem \ref{t4.1} and if%
\begin{equation}
\func{Re}\left\langle \sum_{i\in F}\Phi _{i}e_{i}-\left( \lambda x\pm \left(
1-\lambda \right) y\right) ,\lambda x\pm \left( 1-\lambda \right)
y-\sum_{i\in F}\phi _{i}e_{i}\right\rangle \geq 0  \label{4.8}
\end{equation}
or, equivalently%
\begin{equation}
\left\Vert \lambda x\pm \left( 1-\lambda \right) y-\sum_{i\in F}\frac{\Phi
_{i}+\phi _{i}}{2}\cdot e_{i}\right\Vert \leq \frac{1}{2}\left( \sum_{i\in
F}\left\vert \Phi _{i}-\phi _{i}\right\vert ^{2}\right) ^{\frac{1}{2}},
\label{4.9}
\end{equation}%
then we have the inequality%
\begin{equation}
\left\vert \func{Re}\left[ \left\langle x,y\right\rangle -\sum_{i\in
F}\left\langle x,e_{i}\right\rangle \left\langle e_{i},y\right\rangle \right]
\right\vert \leq \frac{1}{16}\cdot \frac{1}{\lambda \left( 1-\lambda \right) 
}\sum_{i\in F}\left\vert \Phi _{i}-\phi _{i}\right\vert ^{2}.  \label{4.10}
\end{equation}%
The constant $\frac{1}{16}$ is best possible in (\ref{4.10}).
\end{corollary}

\begin{remark}
\label{r4.3}If $H$ is a real inner product space and $m_{i},M_{i}\in \mathbb{%
R}$ with the property%
\begin{equation}
\left\langle \sum_{i\in F}M_{i}e_{i}-\left( \lambda x\pm \left( 1-\lambda
\right) y\right) ,\lambda x\pm \left( 1-\lambda \right) y-\sum_{i\in
F}m_{i}e_{i}\right\rangle \geq 0  \label{4.11}
\end{equation}%
or, equivalently,%
\begin{equation}
\left\Vert \lambda x\pm \left( 1-\lambda \right) y-\sum_{i\in F}\frac{%
M_{i}+m_{i}}{2}\cdot e_{i}\right\Vert \leq \frac{1}{2}\left[ \sum_{i\in
F}\left( M_{i}-m_{i}\right) ^{2}\right] ^{\frac{1}{2}},  \label{4.12}
\end{equation}%
then we have the Gr\"{u}ss type inequality%
\begin{equation}
\left\vert \left\langle x,y\right\rangle -\sum_{i\in F}\left\langle
x,e_{i}\right\rangle \left\langle e_{i},y\right\rangle \right\vert \leq 
\frac{1}{16}\cdot \frac{1}{\lambda \left( 1-\lambda \right) }\sum_{i\in
F}\left( M_{i}-m_{i}\right) ^{2}.  \label{4.13}
\end{equation}
\end{remark}

\section{Integral Inequalities\label{s5}}

Let $\left( \Omega ,\Sigma ,\mu \right) $ be a measure space consisting of a
set $\Omega ,$ a $\sigma -$algebra of parts $\Sigma $ and a countably
additive and positive measure $\mu $ on $\Sigma $ with values in $\mathbb{R}%
\cup \left\{ \infty \right\} .$ Let $\rho \geq 0$ be a $\mu -$measurable
function on $\Omega .$ Denote by $L_{\rho }^{2}\left( \Omega ,\mathbb{K}%
\right) $ the Hilbert space of all real or complex valued functions defined
on $\Omega $ and $2-\rho -$integrable on $\Omega ,$ i.e.,%
\begin{equation}
\int_{\Omega }\rho \left( s\right) \left\vert f\left( s\right) \right\vert
^{2}d\mu \left( s\right) <\infty .  \label{5.1}
\end{equation}

Consider the family $\left\{ f_{i}\right\} _{i\in I}$ of functions in $%
L_{\rho }^{2}\left( \Omega ,\mathbb{K}\right) $ with the properties that%
\begin{equation}
\int_{\Omega }\rho \left( s\right) f_{i}\left( s\right) \overline{f_{j}}%
\left( s\right) d\mu \left( s\right) =\delta _{ij},\ \ \ i,j\in I,
\label{5.2}
\end{equation}%
where $\delta _{ij}$ is $0$ if $i\neq j$ and $\delta _{ij}=1$ if $i=j.$ $%
\left\{ f_{i}\right\} _{i\in I}$ is an orthornormal family in $L_{\rho
}^{2}\left( \Omega ,\mathbb{K}\right) .$

The following proposition holds.

\begin{proposition}
\label{p5.1}Let $\left\{ f_{i}\right\} _{i\in I}$ be an orthornormal family
of functions in $L_{\rho }^{2}\left( \Omega ,\mathbb{K}\right) ,$ $F$ a
finite subset of $I,$ $\phi _{i},\Phi _{i}\in \mathbb{K}$ $\left( i\in
F\right) $ and $f\in L_{\rho }^{2}\left( \Omega ,\mathbb{K}\right) ,$ so
that either%
\begin{equation}
\int_{\Omega }\rho \left( s\right) \func{Re}\left[ \left( \sum_{i\in F}\Phi
_{i}f_{i}\left( s\right) -f\left( s\right) \right) \left( \overline{f}\left(
s\right) -\sum_{i\in F}\overline{\phi _{i}}\text{ }\overline{f_{i}}\left(
s\right) \right) \right] d\mu \left( s\right) \geq 0  \label{5.3}
\end{equation}%
or, equivalently,%
\begin{equation}
\int_{\Omega }\rho \left( s\right) \left\vert f\left( s\right) -\sum_{i\in F}%
\frac{\Phi _{i}+\phi _{i}}{2}\cdot f_{i}\left( s\right) \right\vert ^{2}d\mu
\left( s\right) \leq \frac{1}{4}\sum_{i\in F}\left\vert \Phi _{i}-\phi
_{i}\right\vert ^{2}.  \label{5.4}
\end{equation}%
Then we have the inequality%
\begin{align}
0& \leq \int_{\Omega }\rho \left( s\right) \left\vert f\left( s\right)
\right\vert ^{2}d\mu \left( s\right) -\sum_{i\in F}\left\vert \int_{\Omega
}\rho \left( s\right) f\left( s\right) \overline{f_{i}}\left( s\right) d\mu
\left( s\right) \right\vert ^{2}  \label{5.5} \\
& \leq \frac{1}{4}\sum_{i\in F}\left\vert \Phi _{i}-\phi _{i}\right\vert
^{2}-\sum_{i\in F}\left\vert \frac{\Phi _{i}+\phi _{i}}{2}-\int_{\Omega
}\rho \left( s\right) f\left( s\right) \overline{f_{i}}\left( s\right) d\mu
\left( s\right) \right\vert ^{2}  \notag \\
& \leq \frac{1}{4}\sum_{i\in F}\left\vert \Phi _{i}-\phi _{i}\right\vert
^{2}.  \notag
\end{align}%
The constant $\frac{1}{4}$ is best possible in both inequalities.
\end{proposition}

The proof follows by Theorem \ref{t2.2} applied for the Hilbert space $%
L_{\rho }^{2}\left( \Omega ,\mathbb{K}\right) $ and the orthornormal family $%
\left\{ f_{i}\right\} _{i\in I}.$

The following Gr\"{u}ss type inequality also holds.

\begin{proposition}
\label{p5.2}Let $\left\{ f_{i}\right\} _{i\in I}$ and $F$ be as in
Proposition \ref{p5.1}. If $\phi _{i},\Phi _{i},\gamma _{i},\Gamma _{i}\in 
\mathbb{K}$ $\left( i\in F\right) $ and $f,g\in L_{\rho }^{2}\left( \Omega ,%
\mathbb{K}\right) $ so that either%
\begin{align}
\int_{\Omega }\rho \left( s\right) \func{Re}\left[ \left( \sum_{i\in F}\Phi
_{i}f_{i}\left( s\right) -f\left( s\right) \right) \left( \overline{f}\left(
s\right) -\sum_{i\in F}\overline{\phi _{i}}\text{ }\overline{f_{i}}\left(
s\right) \right) \right] d\mu \left( s\right) & \geq 0,  \label{5.6} \\
\int_{\Omega }\rho \left( s\right) \func{Re}\left[ \left( \sum_{i\in
F}\Gamma _{i}f_{i}\left( s\right) -g\left( s\right) \right) \left( \overline{%
g}\left( s\right) -\sum_{i\in F}\overline{\gamma _{i}}\text{ }\overline{f_{i}%
}\left( s\right) \right) \right] d\mu \left( s\right) & \geq 0,  \notag
\end{align}%
or, equivalently,%
\begin{align}
\int_{\Omega }\rho \left( s\right) \left\vert f\left( s\right) -\sum_{i\in F}%
\frac{\Phi _{i}+\phi _{i}}{2}f_{i}\left( s\right) \right\vert ^{2}d\mu
\left( s\right) & \leq \frac{1}{4}\sum_{i\in F}\left\vert \Phi _{i}-\phi
_{i}\right\vert ^{2},  \label{5.7} \\
\int_{\Omega }\rho \left( s\right) \left\vert g\left( s\right) -\sum_{i\in F}%
\frac{\Gamma _{i}+\gamma _{i}}{2}f_{i}\left( s\right) \right\vert ^{2}d\mu
\left( s\right) & \leq \frac{1}{4}\sum_{i\in F}\left\vert \Gamma _{i}-\gamma
_{i}\right\vert ^{2},  \notag
\end{align}%
then we have the inequalities%
\begin{multline}
\left\vert \int_{\Omega }\rho \left( s\right) f\left( s\right) \overline{%
g\left( s\right) }d\mu \left( s\right) \right.   \label{5.8} \\
-\left. \sum_{i\in F}\int_{\Omega }\rho \left( s\right) f\left( s\right) 
\overline{f_{i}}\left( s\right) d\mu \left( s\right) \int_{\Omega }\rho
\left( s\right) f_{i}\left( s\right) \overline{g\left( s\right) }d\mu \left(
s\right) \right\vert 
\end{multline}%
\begin{align*}
& \leq \frac{1}{4}\left( \sum_{i\in F}\left\vert \Phi _{i}-\phi
_{i}\right\vert ^{2}\right) ^{\frac{1}{2}}\left( \sum_{i\in F}\left\vert
\Gamma _{i}-\gamma _{i}\right\vert ^{2}\right) ^{\frac{1}{2}} \\
& \ \ \ \ \ \ \ \ -\sum_{i\in F}\left\vert \frac{\Phi _{i}+\phi _{i}}{2}%
-\int_{\Omega }\rho \left( s\right) f\left( s\right) \overline{f_{i}}\left(
s\right) d\mu \left( s\right) \right\vert \left\vert \frac{\Gamma
_{i}+\gamma _{i}}{2}-\int_{\Omega }\rho \left( s\right) g\left( s\right) 
\overline{f_{i}}\left( s\right) d\mu \left( s\right) \right\vert  \\
& \leq \frac{1}{4}\left( \sum_{i\in F}\left\vert \Phi _{i}-\phi
_{i}\right\vert ^{2}\right) ^{\frac{1}{2}}\left( \sum_{i\in F}\left\vert
\Gamma _{i}-\gamma _{i}\right\vert ^{2}\right) ^{\frac{1}{2}}.
\end{align*}%
The constant $\frac{1}{4}$ is the best possible.
\end{proposition}

The proof follows by Theorem \ref{t3.1} and we omit the details.

\begin{remark}
\label{r5.3}Similar results may be stated if one applies the inequalities in
Section \ref{s4}. We omit the details.
\end{remark}

In the case of real spaces, the following corollaries provide much simpler
sufficient conditions for the counterpart of Bessel's inequality (\ref{5.5})
or for the Gr\"{u}ss type inequality (\ref{5.8}) to hold.

\begin{corollary}
\label{c5.4}Let $\left\{ f_{i}\right\} _{i\in I}$ be an orthornormal family
of functions in the real Hilbert space $L_{\rho }^{2}\left( \Omega \right) ,$
$F$ a finite part of $I,$ $M_{i},m_{i}\in \mathbb{R}$ \ $\left( i\in
F\right) $ and $f\in L_{\rho }^{2}\left( \Omega \right) $ so that%
\begin{equation}
\sum_{i\in F}m_{i}f_{i}\left( s\right) \leq f\left( s\right) \leq \sum_{i\in
F}M_{i}f_{i}\left( s\right) \text{ \ for \ }\mu -\text{a.e. }s\in \Omega ,
\label{5.9}
\end{equation}%
then we have the inequalities%
\begin{align}
0& \leq \int_{\Omega }\rho \left( s\right) f^{2}\left( s\right) d\mu \left(
s\right) -\sum_{i\in F}\left[ \int_{\Omega }\rho \left( s\right) f\left(
s\right) f_{i}\left( s\right) d\mu \left( s\right) \right] ^{2}  \label{5.10}
\\
& \leq \frac{1}{4}\sum_{i\in F}\left( M_{i}-m_{i}\right) ^{2}-\sum_{i\in F}%
\left[ \frac{M_{i}+m_{i}}{2}-\int_{\Omega }\rho \left( s\right) f\left(
s\right) f_{i}\left( s\right) d\mu \left( s\right) \right] ^{2}  \notag \\
& \leq \frac{1}{4}\sum_{i\in F}\left( M_{i}-m_{i}\right) ^{2}.  \notag
\end{align}%
The constant $\frac{1}{4}$ is best possible.
\end{corollary}

\begin{corollary}
\label{c5.5}Let $\left\{ f_{i}\right\} _{i\in I}$ and $F$ be as in Corollary %
\ref{c5.4}. If $M_{i},m_{i},N_{i},n_{i}\in \mathbb{R}$ $\left( i\in F\right) 
$ and $f,g\in L_{\rho }^{2}\left( \Omega \right) $ such that%
\begin{equation}
\sum_{i\in F}m_{i}f_{i}\left( s\right) \leq f\left( s\right) \leq \sum_{i\in
F}M_{i}f_{i}\left( s\right)  \label{5.11}
\end{equation}%
and%
\begin{equation*}
\sum_{i\in F}n_{i}f_{i}\left( s\right) \leq g\left( s\right) \leq \sum_{i\in
F}N_{i}f_{i}\left( s\right) \text{ \ for \ }\mu -\text{a.e. }s\in \Omega .
\end{equation*}%
Then we have the inequalities%
\begin{align}
& \left\vert \int_{\Omega }\rho \left( s\right) f\left( s\right) g\left(
s\right) d\mu \left( s\right) \right.  \label{5.12} \\
& \ \ \ \ \ \ \ \ \ \ \ \ -\left. \sum_{i\in F}\int_{\Omega }\rho \left(
s\right) f\left( s\right) f_{i}\left( s\right) d\mu \left( s\right)
\int_{\Omega }\rho \left( s\right) g\left( s\right) f_{i}\left( s\right)
d\mu \left( s\right) \right\vert  \notag \\
& \leq \frac{1}{4}\left( \sum_{i\in F}\left( M_{i}-m_{i}\right) ^{2}\right)
^{\frac{1}{2}}\left( \sum_{i\in F}\left( N_{i}-n_{i}\right) ^{2}\right) ^{%
\frac{1}{2}}  \notag \\
& \ \ \ \ \ \ \ \ \ -\sum_{i\in F}\left\vert \frac{M_{i}+m_{i}}{2}%
-\int_{\Omega }\rho \left( s\right) f\left( s\right) f_{i}\left( s\right)
d\mu \left( s\right) \right\vert  \notag \\
& \ \ \ \ \ \ \ \ \ \ \ \ \ \ \ \ \ \ \ \ \ \ \ \ \ \ \ \ \ \times
\left\vert \frac{N_{i}+n_{i}}{2}-\int_{\Omega }\rho \left( s\right) g\left(
s\right) f_{i}\left( s\right) d\mu \left( s\right) \right\vert  \notag \\
& \leq \frac{1}{4}\left( \sum_{i\in F}\left( M_{i}-m_{i}\right) ^{2}\right)
^{\frac{1}{2}}\left( \sum_{i\in F}\left( N_{i}-n_{i}\right) ^{2}\right) ^{%
\frac{1}{2}}.  \notag
\end{align}
\end{corollary}

\end{document}